\documentclass[11pt]{amsart}

\usepackage[square,compress,comma, numbers,sort]{natbib}
\usepackage[colorlinks=true, citecolor=red, linkcolor=blue]{hyperref}
\usepackage{amsfonts,mathtools}

\allowdisplaybreaks[4]

\usepackage{amsmath}
\usepackage{bbm}
\usepackage{amssymb}
\usepackage{mathtools}

\usepackage{color}
\definecolor{c20}{rgb}{0.,0.7,0.}
\definecolor{c30}{rgb}{0.,0.,1.}
\definecolor{c40}{rgb}{1,0.1,0.7}
\definecolor{c50}{rgb}{1,0,0}
\definecolor{c60}{rgb}{1,0.9,0.1}

\def\Var{\text{Var}}

\definecolor{c20}{rgb}{0.,0.7,0.}
\definecolor{c30}{rgb}{0.,0.,1.}
\definecolor{c40}{rgb}{1,0.1,0.7}
\definecolor{c50}{rgb}{1,0,0}
\definecolor{c60}{rgb}{1,0.9,0.1}

\newcommand{\md}{\mathbb{D}_H}
\newcommand{\bp}{\bar{\psi}}
\newcommand{\hu}{\hat{u}}
\newcommand{\ve}{\varepsilon}

\newcommand{\abs}[1]{\left\lvert #1 \right\rvert}

\newcommand{\E}[1]{\mathbb{E}\left\{ #1\right\}}

\newcommand{\pk}[1]{\mathbb{P} \left\{ #1 \right \} }
\newcommand{\PP}{\mathbb{P}}

\newcommand{\R}{\mathbb{R}}

\newcommand{\inr}{\in \R}

\newcommand{\limit}[1]{\lim_{#1 \to   \infty}}

\newcommand{\BQN}{\begin{eqnarray}}
\newcommand{\EQN}{\end{eqnarray}}
\newcommand{\BQNY}{\begin{eqnarray*}}
\newcommand{\EQNY}{\end{eqnarray*}}

\newcommand{\BS}{\begin{sat}}
\newcommand{\ES}{\end{sat}}
\newcommand{\BT}{\begin{theo}}
\newcommand{\ET}{\end{theo}}
\newcommand{\BK}{\begin{korr}}
\newcommand{\EK}{\end{korr}}

\newcommand{\PH}{\overline{\Phi}}
\DeclareMathOperator{\cov}{cov}

\newcommand{\BD}{\begin{de}}
\newcommand{\ED}{\end{de}}

\newcommand{\BIT}{\begin{itemize}}
\newcommand{\EIT}{\end{itemize}}
\newcommand{\BDI}{\begin{description}}
\newcommand{\EDI}{\end{description}}

\newcommand{\BRM}{\begin{remarks}}
\newcommand{\ERM}{\end{remarks}}

\newcommand{\BEL}{\begin{lem}}
\newcommand{\EEL}{\end{lem}}

\def\bqny#1{{\begin{eqnarray*} #1 \end{eqnarray*}}}
\def\bqn#1{{\begin{eqnarray} #1 \end{eqnarray}}}

\newtheorem{theo}{Theorem}[section]
\newtheorem{sat}[theo]{Proposition}
\newtheorem{de}[theo]{Definition}
\newtheorem{lem}[theo]{Lemma}

\newtheorem{korr}[theo]{Corollary}
\newtheorem{remark}[theo]{Remark}
\newtheorem{remarks}[theo]{Remarks}

\newcommand{\COM}[1]{}

\newcommand{\h}{ \mathcal{H} }

\newcommand{\ttt}{\hat{t}}
\newcommand{\sss}{\hat{s}}

\def\IF{\infty}

\newcommand{\QED}{\hfill $\Box$}

\def\ve{\varepsilon}

\def\IF{\infty}

\begin{document}

\title[]{Simultaneous Ruin Probability for Two-Dimensional Fractional Brownian motion risk Process over Discrete Grid}

\author[]{Grigori Jasnovidov}
\address{Grigori Jasnovidov, Department of Actuarial Science, 
University of Lausanne,\\
UNIL-Dorigny, 1015 Lausanne, Switzerland
}
\email{Grigori.Jasnovidov@unil.ch}

\bigskip

\date{\today}
 \maketitle

{\bf Abstract:} This paper derives the asymptotic behavior of
the following ruin probability
$$\pk{\exists t \in G(\delta):B_H(t)-c_1t>q_1u,B_H(t)-c_2t>q_2u},
 \ \ \ u \to \IF,$$
where $B_H$ is a standard fractional Brownian motion, $c_1,q_1,c_2,q_2>0$
and $G(\delta)$ denotes a regular grid
$\{0,\delta, 2\delta,...\}$ for some $\delta>0$.
The approximation depends on $H$, $\delta$ (only when $H\leq 1/2$) and the relations between parameters $c_1,q_1,c_2,q_2$.

\section{Introduction}
Let $B_H(t), \ t \inr $ be a standard fractional Brownian motion (fBM)
with zero-mean and covariance function
$$\cov(B_H(t),B_H(s))=\frac{1}{2}(|t|^{2H}+|s|^{2H}-|t-s|^{2H}), \ \ \ \
H \in (0,1), \ \ \ s,t \inr.$$
Define two risk processes
$$
R^{(H)}_{i,u}(t) = q_iu+c_it-B_H(t), \ \ i=1,2,
$$
where $c_i,q_i>0$.
The discrete simultaneous ruin probability over the infinite time
horizon is defined by
\bqn{\label{prob}
\bp_{\delta,H}(u) = \pk{\exists t \in G(\delta):R^{(H)}_{1,u}(t)<0,R^{(H)}_{2,u}(t)<0},
}
where $G(\delta)$ denotes the grid $\{0,\delta, 2\delta,...\}$
(if $\delta=0$, then set $G(\delta)=[0,\IF)$). For positive $\delta$
the simultaneous ruin probability is of interest both for theory-oriented studies and for
applications in reinsurance (see e.g. \cite{Lanpeng2BM} and references therein).
In this paper we investigate only the discrete setup; the continuous problem
has been already solved in \cite{Lanpeng2BM}.
For any possible choices of positive $\delta$ and $H\in (0,1)$
it is not possible to calculate $\bp_{\delta,H}(u)$
explicitly. A natural question when lack of explicit
formulas is the case, is how can we approximate $\bp_{\delta,H}(u)$
for large $u$? Also of interest is to know what is the role of $\delta$,
does it affects the ruin probability in the considered risk model?
Theorem \ref{maintheo} gives detailed answers for these questions.
Our results show that the discrete time ruin probabilities
behave differently from continuous if $H\le 1/2$.
\\

\COM{
To prove our results we mainly use the Double-Sum method. We consider different scenarios, which are determined by value of $H$ (namely, $H<1/2$ or $H=1/2$ or $H>1/2$) and the relations between $c_1,q_1,c_2,q_2$. When $H=1/2$ we use the self-similarity, stationarity and independence of the increments of Brownian motion (BM), however
when $H\neq 1/2$ we utilize only the self-similarity property.\\
}

Also of certain interest is the finite time horizon setup of the problem.
For fixed $T>0$
the discrete simultaneous ruin probability over a finite time horizon is
\bqn{\label{vscontin}
\bar{\zeta}_{H,T}(u) = \pk{\exists t \in [0,T]:R^{(H)}_{1,u}(t)<0,R^{(H)}_{2,u}(t)<0}
.}

The corresponding discrete ruin problem over a finite time horizon
(in contrary with the infinite case) is trivial, since set
$[0,T]\cap G(\delta)$ consists of finite number of elements and hence asymptotics
of the large deviation is determined by the unique maximizer
of the variance of the
process (this e.g., follows immediately from Lemma 2.3
in \cite{PickandsA} or Proposition 2.4.2 in \cite{20lectures}). In the
light of this we shall be concerned only with the continuous time problem.
Asymptotics of $\bar{\zeta}_{H,T}(u)$ are discussed in Remark \ref{prop}.
\\
\\
We organize the paper in the following way.
The next section gives notation, necessary assumptions and the main results.
All proofs are relegated to Section 3, while some technical calculations
are presented in the Appendix.

\section{Main Results}

First of all we eliminate the  trends via self-similarity of
fBM. We have for any $u>0$
\bqny{\bar{\psi}_{\delta,H}(u)
&=&\pk{\exists t \in G(\delta): B_H(t)>q_1u+c_1t,B_H(t)>q_2u+c_2t}
\\ &=&
\pk{\exists tu \in G(\delta):B_H(tu)>(q_1+c_1t)u,B_H(tu)>(q_2+c_2t)u}
\\&=&
\pk{\exists t \in G(\delta/u): \frac{B_H(t)}{c_1t+q_1}>u^{1-H},\frac{B_H(t)}{c_2t+q_2}>u^{1-H}}
\\ &=&
\pk{\exists t \in G(\delta/u): \frac{B_H(t)}{\max(c_1t+q_1,c_2t+q_2)}>u^{1-H}}.
}
If the two lines $q_1 + c_1t$ and $q_2 + c_2t$ do
not intersect over $(0,\infty)$, then the problem degenerates to
the one-dimensional case, which is discussed in Theorem \ref{H1}.
In consideration of that dealing with $\bar{\psi}_{\delta,H}(u)$
we always suppose that
\bqn{\label{cq}
c_1>c_2, \ \ q_2> q_1.
}
It turns out, that the variance of
$\frac{B_H(t)}{\max(c_1t+q_1,c_2t+q_2)}$
can achieve its unique maxima only at one of the following points:
\bqn{\label{t1t2}
t_1 = \frac{Hq_1}{c_1(1-H)}, \ \ \ \ t_2 =\frac{Hq_2}{c_2(1-H)}, \ \ \ \ t^* = \frac{q_2-q_1}{c_1-c_2}
.}
It follows from \eqref{cq} that
$t_1<t_2.$
As we show later, the order between $t_1,t_2$ and $t^*$ determines
the asymptotics of $\bp_{\delta,H}(u)$ as $u \to \IF$.
\\
\\
Denote by $\Phi$ and $\PH$ the distribution and survival functions of a
standard normal random variable, respectively.
 For notational simplicity we shall write
$B(t)$ instead of $B_{1/2}(t)$ and $\bar\psi_\delta(u)$
instead of $\bar{\psi}_{\delta,1/2}(u)$.\\
Define Pickands constant $\mathbb{H}_{2H}$ by
$$\mathbb{H}_{2H} = \limit{S}\frac{1}{S}
\E{\sup\limits_{t \in [0,S]}e^{\sqrt 2B_H(t)-|t|^{2H}}}, \ \ \ \ \ \
H\in(0,1)$$
and for any $\alpha >0$ the discrete Pickands constant by
$$\mathcal{H}_\alpha=	\E{ \frac{\sup_{t\in \alpha
\mathbb{Z}} e^{ \sqrt{2} B(t)- \abs{t}} }{\alpha \sum\limits_
{t \in \alpha \mathbb{Z}}e^{ \sqrt{2} B(t)- \abs{t}}  } }.$$
It is known that $\mathbb{H}_{2H}$ and $\mathcal{H}_\alpha$ are positive and finite constants (see \cite{PickandsA}, \cite{20lectures}, \cite{DiekerY}).
For any real numbers $a\le b$ and $\alpha>0$ denote $[a,b]_\alpha = [a,b] \bigcap \alpha \mathbb{Z}$.
Define
$$\h_\alpha^k[S,T] = \E{
\sup\limits_{t \in [S,T]_\alpha} e^{\sqrt 2 B(t)-|t|+k(t)}}$$
for some function $k(t)$ and $-\IF <S< T< \IF$. Denote further
$$\h_\alpha^k = \lim\limits_{T \rightarrow \IF}\h_\alpha^k[-T,T]
$$
when the limit exists. We refer to \cite{20lectures} for properties
of Piterbarg constants.
Let
$$d(t)  =
\mathbb{I}(t<0)\frac{(q_2c_1+c_2q_1-2q_2c_2)t}{c_1q_2-q_1c_2} +
\mathbb{I}(t\geq 0)\frac{(2c_1q_1-c_1q_2-q_1c_2)t}{c_1q_2-q_1c_2}$$
and
$$  d_\delta(t) =
\mathbb{I}(t<0)\frac{(q_2c_1+c_2q_1-2q_2c_2)t}{c_1q_2-q_1c_2}  +
\mathbb{I}(t\geq 0)\big(\frac{(2c_1q_1-c_1q_2-q_1c_2)t}{c_1q_2-q_1c_2}-
\delta\frac{(c_1q_2-q_1c_2)(c_1-c_2)}{q_2-q_1}\big),
$$
where $\mathbb{I}(\cdot)$ is the indicator function.
Define constants
$$C^{(i)}_H = \frac{c_i^Hq_i^{1-H}}{H^H(1-H)^{1-H}}, \ \ \ i =1,2.$$
The theorem below establishes the asymptotics of
$\bar{\psi}_{\delta,H}(u)$.

\begin{theo} \label{maintheo} Let $\delta>0$ and
$u \rightarrow \IF$. \\
1) If $t^*\notin (t_1,t_2)$, then
\bqny{
\bar{\psi}_{\delta,H}(u) \sim
\begin{cases}
\mathbb{H}_{2H}\frac{2^{\frac{1}{2}-\frac{1}{2H}}\sqrt \pi}{H^{1/2}(1-H)^{1/2}}(C_{H}^{(i)}u^{1-H})^{\frac{1}{H}-1}
\PH(C_{H}^{(i)}u^{1-H}), & H>1/2\\
\mathcal{H}_{2c_i^2\delta}e^{-2c_iq_iu}, & H=1/2 \\
\frac{\sqrt{2\pi} H^{H+1/2} q_i^H u^H}{\delta c_i^{H+1}(1-H)^{H+1/2}}
\PH(C^{(i)}_Hu^{1-H}), & H<1/2,
\end{cases}
}
where $i=1$ if $t^*<t_1$ and $i=2$ if $t^*>t_2$.\\
2) If $t^*\in (t_1,t_2)$, then
with $\md = \frac{c_1t^*+q_1}{(t^*)^H} = \frac{c_2t^*+q_2}{(t^*)^H}$
when $H>1/2$
\bqn{\label{claimHsim}
\bar{\psi}_{\delta,H}(u) \sim \PH(\md u^{1-H})
,}
when $H=1/2$
\bqn{\label{claim_case2_H=1/2_maintheo}
\mathcal{H}^{d_\delta}_\gamma
 \PH(\mathbb{D}_{1/2}\sqrt u)(1+o(1))
  \leq
\bp_\delta(u) \leq
 A\mathcal{H}^d_\gamma
\PH(\mathbb{D}_{1/2}\sqrt u)(1+o(1)),}

where $\mathcal{H}^d_\gamma,\mathcal{H}^{d_\delta}_\gamma \in (0,\IF)$ and
\bqn{\label{eta,A,gamma_definition}
A  =
e^{\delta\frac{(c_1q_2-c_2q_1)(c_1q_2+q_1c_2-2c_2q_2)}{2(q_2-q_1)^2}}, \ \ \ \
 \ \ \gamma = \frac{\delta(c_1q_2-q_1c_2)^2}{2(q_2-q_1)^2},
}
when $H<1/2$
\bqn{ \label{small}
2e^{-Bu^{1-H}}\PH(\md u^{1-H})
(1+o(1))
\leq \bar{\psi}_{\delta,H}(u)
\leq
\PH(\md u^{1-H})(1+o(1))
,}
where
\bqn{\label{ww}
w_1(t) = \frac{(q_1+c_1t)^2}{t^{2H}}, \ \ \ \
w_2(t) = \frac{(q_2+c_2t)^2}{t^{2H}}, \ \ \ \
B = -\frac{\delta w_1'(t^*)w_2'(t^*)}{2(w_1'(t^*)-w_2'(t^*))}>0.}
3) If $t^*=t_i, \ i=1,2$, then
\bqny{
\bar{\psi}_{\delta,H}(u) \sim \frac{1}{2}\times
\begin{cases}
\mathbb{H}_{2H}\frac{2^{\frac{1}{2}-\frac{1}{2H}}\sqrt \pi}{H^{1/2}(1-H)^{1/2}}(C_{H}^{(i)}u^{1-H})^{\frac{1}{H}-1}
\PH(C_{H}^{(i)}u^{1-H}), & H>1/2\\
\mathcal{H}_{2c_i^2\delta}e^{-2c_iq_iu}, & H=1/2 \\
\frac{\sqrt{2\pi} H^{H+1/2} q_i^Hu^H}
{\delta c_i^{H+1}(1-H)^{H+1/2}}
\PH(C^{(i)}_Hu^{1-H}), & H<1/2.
\end{cases}
}

\end{theo}

\begin{remark}\label{rem2}
The bounds in \eqref{small}
are exact. Namely, there exist two tending to infinity sequences
$\{u_n\}_{n\in \mathbb{N}}$ and $\{v_n\}_{n\in \mathbb{N}}$
such that as $n\to \IF$
\bqny{
\bar{\psi}_{\delta,H}(u_n) \sim \PH(\md u_n^{1-H}), \quad
\bar{\psi}_{\delta,H}(v_n) \sim
2e^{-B v_n^{1-H}}\PH(\md v_n^{1-H}).
}
\end{remark}

\begin{remark}
Let $u \to \IF$. Then for case 1) it holds
$$\bar{\psi}_{\delta,H}(u) \sim
\pk{\sup\limits_{t \in G(\delta)}(B_H(t)-c_it)>q_iu}  \ \ \ i = 1,2,$$
where $i=1$ if $t^*<t_1$ and $i=2$ if $t^*>t_2$.
For case 3) it holds
$$\bar{\psi}_{\delta,H}(u) \sim  \frac{1}{2}\pk{\sup\limits_{t \in G(\delta)}(B_H(t)-c_it)>q_iu} \ \ \ i = 1,2.$$
\end{remark}
To study the asymptotics of the two-dimensional ruin probability over
the infinite time horizon
crucial is the asymptotic approximation of the one-dimensional one.
The following auxiliary theorem derives the approximations of this
probability.
\BT \label{H1}
For any $\delta>0$ with $C_H = \frac{c^H}{H^H(1-H)^{1-H}}$
as $u \to \IF$
\bqn{\label{one-dimensional_theorem} \ \ \ \
\pk{\exists t \in G(\delta): B_H(t)-ct>u} \sim
\begin{cases}
\mathbb{H}_{2H}\frac{2^{\frac{1}{2}-\frac{1}{2H}}\sqrt \pi}
{H^{1/2}(1-H)^{1/2}}
(C_Hu^{1-H})^{1/H-1}\PH(C_Hu^{1-H}), & H>1/2,\\
\mathcal{H}_{2c^2\delta}  e^{-2cu}, & H=1/2\\
\frac{\sqrt{2\pi} H^{H+1/2} u^H}{\delta c^{H+1}(1-H)^{H+1/2}}
\PH(C_Hu^{1-H}), & H<1/2.
\end{cases}
}
\ET

\begin{remark}\label{remcompasymp}
When $H>1/2$ the asymptotics of
the discrete probabilities in Theorems \ref{maintheo} and \ref{H1}
are the same as in the continuous case and do not depend on $\delta$.
The cause is "high density" of the grid $G(\delta/u)$.
If $H=1/2$ the asymptotics differ only in the constants.
For $H<1/2$ the discrete asymptotics are infinitely smaller
then the corresponding continuous.
Here the grid $G(\delta/u)$ is "sparse".
All these statements directly follow from Proposition 2.1
and Theorem 3.1 in \cite{Lanpeng2BM}.
\end{remark}

Next we study the finite-horizon case.
Here for large $u$
the two-dimensional ruin probability always
reduces to the one-dimensional one, that has been already studied in
\cite{SumFBMDebicki},\cite{DebRolski2002}. More precisely, we have

\begin{remark}\label{prop}
Regardless of \eqref{cq} for any $T>0$
with $\lambda(u)=\frac{\max(q_1u+c_1T,q_2u+c_2T)}{T^{H}}$
as $u \to \IF$ \\
\bqny{
\bar{\zeta}_{H,T}(u) \sim
\begin{cases}
\mathbb{H}_{2H}(\lambda(u))^{\frac{1-2H}{H}}\frac{(1/2)^{(1/2H)}}{H}
\PH(\lambda(u)), & H<1/2\\
\PH(\lambda(u)), & H>1/2
\end{cases}
}
and
$$\bar{\zeta}_{\frac{1}{2},T}(u) =
\PH(\frac{uq_i}{\sqrt T}+c_i\sqrt T)+e^{-2c_iq_iu}
\PH(\frac{uq_i}{\sqrt T}-c_i\sqrt T), \quad i =1,2,$$
where $i=1$ if $(q_1,c_1)\ge (q_2,c_2)$ in the alphabetical order and
$i=2$ otherwise.
\end{remark}

\section{Proofs}
Let $\mathcal{N}$ be a standard normal random variable independent of all stochastic processes which we consider. We write $\xi \sim \mathcal{N}(\mu,\sigma^2)$ if $\xi$ is a Gaussian random variable with expectation $\mu$ and variance $\sigma^2$.
We reserve notation $C, C_1$ for some positive constants that do not depend on $u$,
they might differ in different places.
\\
\\
\textbf{Proof of Theorem \ref{maintheo}.}
\textbf{Case (1).} Assume, that $t_1<t^*$, case $t_2>t^*$ follows
by the same arguments. Notice, that
\bqn{\label{psi1def} \ \ \ \ \ \ \ \
\bar{\psi}_{\delta,H}(u) \leq
\pk{\exists t \in G(\delta):B_H(t)-c_1t>q_1u}
&=&\pk{\exists t \in G(\delta/u):\frac{B_H(t)}{c_1t+q_1}>u^{1-H}}
\notag\\&=:&
\psi^{(1)}_{\delta,H}(u).
}
Denote
\bqn{\label{videf}
V_i(t) = \frac{B_H(t)}{c_it+q_i}, \ \ \ i=1,2.
}
Since $t^*<t_1$ we have for any positive $\ve<t_1-t^*$
\bqny{
\bar{\psi}_{\delta,H}(u) &\geq&
\pk{\exists t \in [t_1-\ve,t_1+\ve]_{\frac{\delta}{u}}: V_1(t)>u^{1-H},
V(t_2)>u^{1-H} }
\\ &=&
\pk{\exists t \in [t_1-\ve,t_1+\ve]_{\frac{\delta}{u}}: V_1(t)>u^{1-H}}
\\ &\sim&
\psi^{(1)}_{\delta,H}(u), \ \ \ \ u \to \IF.
}
The last line above follows from
\bqn{\label{borH}
\pk{\exists t \notin [t_1-\ve,t_1+\ve]: V_1(t)>u^{1-H}}
= o(\psi^{(1)}_{\delta,H}(u)), \ \ \ \ u \to \IF,}
proof is in the Appendix. Thus
$$\bar{\psi}_{\delta,H}(u) \sim \psi^{(1)}_{\delta,H}(u), \ \ \ \ \ u \to \IF$$
and by Theorem \ref{H1} the claim is established.
\\
\\
\textbf{Case (2).}
Denote
\bqn{\label{ZZ}
Z_H(t) = \frac{B_H(t)}{\max(c_1t+q_1,c_2t+q_2)} \ \ \ \ \
\text{and} \ \ \ \ \
\sigma^2_H(t) = \Var \{ Z_H(t)\}.
}
Notice, that if $t_1\leq t^*\leq t_2$, then
$t^*$ is the unique maximizer of $\sigma_H(t)$.
Moreover, $\sigma_H(t)$ increases over
$(0,t^*)$ and decreases over $(t^*,\IF)$.
Define
\bqn{\label{tu}
t_u = t^*+\frac{\theta_u}{u} \in G(\frac{\delta}{u}), \ \text{where} \ \theta_u \in [0,\delta).
}
\emph{Case $H>1/2$.}
From Theorem 3.1 case (3), $H>\frac{1}{2}$ in \cite{Lanpeng2BM} it follows
that
\bqn{\label{star}
\bar{\psi}_{\delta,H}(u) \leq  \PH(\md u^{1-H})(1+o(1)), \ \ \ \ u \to \IF.
}
We have by the asymptotic ratio $\PH(x) \sim \frac{1}{\sqrt{2\pi}x}e^{-x^2/2}, \ x \to \IF$ (see e.g. \cite{PickandsA}, Lemma 2.1)
\bqny{
\bar{\psi}_{\delta,H}(u) = \pk{\exists t \in G(\frac{\delta}{u}): Z_H(t)>u^{1-H}}
&\geq&
\pk{Z_H(t_u)>u^{1-H}}
\\&=&
\PH(u^{1-H}\frac{q_1+c_1(t^*+\theta_u/u)}{(t^*+\theta_u/u)^H})
\sim \PH(\md u^{1-H}), \ \ u \to \IF.
}

Combining the statement above with \eqref{star} we establish the claim.
\\
\\
\emph{Case $H=1/2$.}
For notation simplicity, here we write $Z(t)$ instead of $Z_{1/2}(t)$.
We have (proof is in the Appendix)
\bqn{\label{informinterv}
\bp_\delta(u) \sim \pk{\exists t \in [t^*-\frac{\ln u}{\sqrt u}, t^*+\frac{\ln u}{\sqrt u}]_{\frac{\delta}{u}} : Z(t)> \sqrt u}, \ \ u \to \IF.}

Next for any
fixed $S,u>0$ we consider the intervals
\BQNY
\Delta_{j,S,u}= [t_u+j Su^{-1}, t_u+(j+1) Su^{-1}]_\frac{\delta}{u},
 \ \ \ \ -N_u  \leq j \leq N_u,
\EQNY
where $N_u=\lfloor S^{-1} \ln(u) \sqrt{u}\rfloor$
and $\lfloor\cdot\rfloor$ is the ceiling function.
Let
\bqny{p_{j,S,u} = \pk{\sup\limits_{t \in \Delta_{j,S,u}}\frac{B(t)}{c_1t+q_1}>\sqrt u} \text{ for } j \geq 0, \ \ \ \
p_{j,S,u} = \pk{\sup\limits_{t \in \Delta_{j,S,u}}\frac{B(t)}{c_2t+q_2}>\sqrt u} \text{ for } j < 0.}

Denote $\Delta = \Delta_{-1}\bigcup\Delta_0$.
We have
\bqn{\label{sums}
\pk{\sup\limits_{t \in \Delta} Z(t)> \sqrt u}
&\leq&
\pk{\exists t \in [t^*-\frac{\ln u}{\sqrt u}, t^*+\frac{\ln u}{\sqrt u}]_{\frac{\delta}{u}} : Z(t)> \sqrt u}
\notag \\&\leq&
\sum\limits_{j=-N_u}^{-2}p_{j,S,u} +
\sum\limits_{j=1}^{N_u}p_{j,S,u}
+ \pk{\sup \limits_{t \in \Delta} Z(t)> \sqrt u}
.}
We shall compute the asymptotic of each summand in the right part
of the inequality above and then compare the asymptotics.
\\
\\
\emph{Approximation of} $\sum\limits_{j=1}^{N_u}p_{j,S,u}$.
We have
\bqny{
p_{j,S,u} = \pk{\sup\limits_{t \in \Delta_{j,S,u}}\frac{B(t)}{c_1t+q_1}>\sqrt u} =
\pk{\sup\limits_{t \in \Delta_{j,S,u}} (B(t)-\sqrt \alpha \mu t) >\sqrt \alpha}
,}
where
$$\sqrt \alpha = q_1\sqrt u \ \text{ and } \ \mu = \frac{c_1}{q_1}.$$
By independence of the increments of BM with
$\alpha=v^2$, $c_{j,S,v}= t_u+ j S v^{-2}$,
$\hat S = Sq_1^2$, $\hat \delta = \delta q_1^2$
 and $\varphi_{v,j}$ being the
probability density function of $\sqrt{c_{j,S,v}} \mathcal{N}$ we have
\bqny{
	p_{j,S,u}
&=&
\pk{\exists_{t\in\Delta_{j,S,u} }: ( B(t) -  \sqrt{\alpha} \mu t)>
		\sqrt{\alpha}  }
\\&=&
\pk{\exists t \in [t_u+\frac{jS}{u},t_u+\frac{(j+1)S}{u}]_{\frac{\delta}{u}}
:B(t)-B(c_{j,S,v})-\sqrt \alpha \mu(t-c_{j,S,v})+B(c_{j,S,v})-
\sqrt \alpha \mu c_{j,S,v}>\sqrt \alpha
}
\\&=&
\int\limits_\R
\pk{\exists t \in [0,\frac{S}{u}]_{\frac{\delta}{u}}:
B(t)-\sqrt \alpha\mu t-\sqrt \alpha \mu c_{j,S,v}
>\sqrt \alpha- x\lvert \sqrt{ c_{j,S,v}} \mathcal{N} = x}
\varphi_{v,j}(x)dx
\\&=&
\int\limits_\R \pk{\exists tv^2 \in [0,Sq_1^2]_{q_1^2\delta}:
\frac{B(tv^2)}{ v}
-\frac{\sqrt \alpha\mu tv^2}{v^2}
-\sqrt \alpha \mu c_{j,S,v}>\sqrt \alpha- x
}\varphi_{v,j}(x)dx
\\&=& \int_{\R} \pk{\exists t \in [0,\hat S]_{\hat \delta}
:( B(t)/v -  v \mu (c_{j,S,v}+ t/v^2)  > v-x
	} \varphi_{v,j}(x) dx\\
	&=&\frac{1}{v} \int_{\R} \pk{\exists t \in [0,\hat S]_{\hat \delta} :( B(t)/v -  v \mu (c_{j,S,v}+ t/v^2)  > v-(v-x/v)
} \varphi_{v,j}(v-x/v) dx\\
	&=&\frac{1}{v} \int_{\R} \pk{\exists t \in [0,\hat S]_{\hat \delta} :( B(t) -  \mu t)   >x+\mu c_{j,S,v}v^2
} \varphi_{v,j}(v-x/v) dx\\
	&=&\frac{1}{v} \int_{\R} \pk{\exists t \in [0,\hat S]_{\hat \delta} :( B(t) -  \mu t)   >x
} \varphi_{v,j}(v(1+\mu c_{j,S,v}) -  x/v) dx\\
	&=&\frac{ e^{ - v^2(1+\mu c_{j,S,v})^2/(2 c_{j,S,v})} }{v \sqrt{2 \pi c_{j,S,v} }}
\int_{\R} \pk{\exists t \in [0,\hat S]_{\hat \delta}  :( B(t) -  \mu t)   >x
} e^{   x(1+\mu c_{j,S,v}) / c_{j,S,v  }-    x^2/(2 c_{j,S,v} v^2)} dx.
}
By Borell-TIS inequality (Lemma 5.3 in \cite{Lanpeng2BM}),
(see also proof of Theorem 1.1 in \cite{HashorvaGrigori} and proof of
\eqref{integral})
we have with $\chi = \frac{1+\mu t^*}{t^*}$ as $u \to \IF$
\bqny{
& \ &
\int\limits_\R \pk{\exists_{t \in [0,\hat S]_{\hat \delta}}  :( B(t) -  \mu t)   >x
} e^{   x(1+\mu c_{j,S,v}) / c_{j,S,v  }-    x^2/(2 c_{j,S,v} v^2)} dx
\\ &\sim&
\int\limits_\R \pk{\exists t \in [0,\hat S]_{\hat \delta}  :( B(t) -  \mu t)   >x
}e^{\chi x}dx
=: J(S).
}

Clearly, $J(S)$ is a non-decreasing function and by
the explicit formula (see \cite{DeM15})
\bqn{\label{brownian_motion_ruin_prob}
\pk{\exists t\ge 0:B(t)-ct>x}=e^{-2cx}, \ \ \ c,x>0
}

we have
\bqny{
J(S) \le
\int\limits_{-\IF}^0e^{\chi x}dx+
\int\limits_0^\IF
 \pk{\exists t \ge 0 : B(t) -  \mu t   >x
}e^{\chi x}dx
=
\frac{1}{\chi}+\int\limits_0^\IF e^{(-2\mu+\chi)x} dx
<\IF,
}
provided by $\frac{\chi}{\mu}<2$ that follows
from $t_1<t^*$. Thus we have that
\bqn{\label{finitness_J(S)}
\limit{S}J(S) \in (0,\IF).
}
Hence we have as $u \to \IF$ and then $S \to \IF$
\bqny{
\sum\limits_{j=1}^{N_u}p_{j,S,u}
\le
C\frac{1}{v}\sum\limits_{j=1}^{N_u}e^{\frac{-v^2(1+\mu c_{j,S,v})^2}{2c_{j,S,v}}}
= \frac{C}{v }
e^{\frac{-v^2(1+\mu t^*)^2}{2t^*}}
\sum\limits_{j=1}^{N_u} e^{-v^2(\frac{(1+\mu c_{j,S,v})^2}{2c_{j,S,v}}-\frac{(1+\mu t^*)^2}{2t^*})}
.}
Setting
$$a(t)= (1+ \mu t)^2/2t = 1/(2t)+ \mu + \mu^2t/2, \quad a'(t)= (-1/t^2 +\mu^2)/2,$$
we have an expansion
$a(t^*+x)= a(t^*)+xa'(t^*)+O(x^2)$ as $x\rightarrow 0$.
Hence (proof is in the Appendix)
\bqn{\label{derivate}
\sum\limits_{j=1}^{N_u} e^{-v^2(\frac{(1+\mu c_{j,S,v})^2}{2c_{j,S,v}}
-\frac{(1+\mu t^*)^2}{2t^*})}
\sim  \sum\limits_{j=1}^{N_u} e^{-v^2a'(t^*)
(\frac{\theta_u+jS}{u})}, \ \ \ u \to \IF
.}
We have with $\omega = a'(t^*)q_1^2>0$
\bqny{\sum\limits_{j=1}^{N_u} e^{-v^2a'(t^*)(\frac{\theta_u+jS}{u})} =
e^{-\omega\theta_u }  \sum\limits_{j=1}^{N_u}e^{-jS\omega}  \le
Ce^{-\omega S},
\ \ \ S \rightarrow \IF.
}
In the light of the calculations above, we have
\bqn{\label{asympsum2}
\sum\limits_{j=1}^{N_u}p_{j,S,u} \le
C\frac{1}{v }
e^{\frac{-v^2(1+\mu t^*)^2}{2t^*}}e^{-\omega S}
\le C \PH(\mathbb{D}_{1/2}\sqrt u)e^{-\omega S}
.
}

\emph{Approximation of} $\pk{\sup \limits_{t \in \Delta} Z(t)> \sqrt u}$.
Let $B^*(t)$ be an independent copy of BM,
$\phi_u(x)$ be the probability density function
of $\sqrt{ut_u}\mathcal{N}$ and
define
\bqn{\label{eta}
\eta = q_1+c_1t^*=q_2+c_2t^*= \frac{c_1q_2-q_1c_2}{c_1-c_2}
.}
By the self-similarity and independence of the increments of BM
we have as $u \to \IF$
\bqn{\label{gausint}
& \ &
\pk{\sup \limits_{t \in \Delta} Z(t)> \sqrt u}
\notag\\&=&
\pk{ \exists \ttt \in [ut_u-S,ut_u)_\delta : B(\ttt)>q_2u+c_2\ttt
\text{ or } \exists t \in [ut_u,ut_u+S]_\delta : B(t)>q_1u+c_1t}
\notag\\&=&
\PP\{ \exists \ttt \in [ut_u-S,ut_u)_\delta : B(\ttt)>q_2u+c_2\ttt
\notag\\& \ &
\text{ or }  \exists t \in [ut_u,ut_u+S]_\delta : (B(t)-B(ut_u))+B(ut_u)>q_1u+c_1t
\}
\notag\\&=&
\PP\{
\exists \ttt \in [ut_u-S,ut_u)_\delta : B(\ttt)>q_2u+c_2ut_u+c_2(\ttt-ut_u)
\notag\\& \ &
\text{ or } \exists t \in [ut_u,ut_u+S]_\delta :
 B^*(t-ut_u)+B(ut_u)>q_1u+c_1ut_u+c_1(t-ut_u)\}
\notag\\&=&
\int\limits_\R \phi_u(\eta u -x)\times \PP\{
\exists \ttt \in [ut_u-S,ut_u)_\delta : B(\ttt)>q_2u+c_2ut_u+c_2(\ttt-ut_u)
\notag\\& \ &
 \ \ \ \ \ \text{ or }
\exists t \in [ut_u,ut_u+S]_\delta : B^*(t-ut_u)+\eta u -x>
q_1u+c_1ut_u+c_1(t-ut_u)|B(ut_u) = \eta u-x\}dx
\\ &=&
\int\limits_\R\PP\{
\exists \ttt \in [-S,0)_\delta : Z_u(\ttt)-c_2\ttt>x+c_2\theta_u
\text{ or }
\exists t \in [0,S]_\delta : B^*(t)-c_1t>x+c_1\theta_u
\}\phi_u(\eta u -x)dx
\notag
\\ &=&
\frac{e^{\frac{-\eta^2u}{2t_u}}}{\sqrt{2\pi ut_u}}
\int\limits_\R\pk{
\exists \ttt \in [-S,0)_\delta : Z_u(\ttt)-c_2\ttt>x+c_2\theta_u \text{ or }
\exists t \in [0,S]_\delta : B^*(t)-c_1t>x+c_1\theta_u
}e^{\frac{\eta x}{t_u}-\frac{x^2}{2ut_u}}dx
\notag\\&\sim&
\frac{1}{\sqrt{2\pi ut^*}}e^{\frac{-\eta^2u}{2t^*}}
e^{\frac{\eta^2\theta_u}{2(t^*)^2}-\frac{\eta c_2\theta_u}{t^*}}
\notag\\&\times&
\int\limits_\R\pk{
\exists \ttt \in [-S,0)_\delta : Z_u(\ttt)-c_2\ttt>x \text{ or }
\exists t \in [0,S]_\delta : B^*(t)-c_1t>x+(c_1-c_2)\theta_u
}e^{\frac{\eta x}{t_u}-\frac{(x-c_2\theta_u)^2}{2ut_u}}dx,
\notag}
where $Z_u(\ttt)$ is a Gaussian process with expectation, covariance and variance defined below ($\sss \leq \ttt$):
$$\E{Z_u(\ttt)} = \frac{\eta u-x}{ut_u}\ttt, \ \ \ \ \ \ \Var \{Z_u(\ttt)\} = -\ttt-\frac{\ttt^2}{ut_u}, \ \ \ \ \ \ \cov (Z_u(\sss),Z_u(\ttt)) = \frac{-\sss\ttt}{ut_u}-\ttt.$$
Technical details of \eqref{gausint} are in the Appendix.
 Since $\eta-2t^*c_2>0$ we have
\bqn{ \label{sandbounds} & \ &
 \ \ \ \ \ \ \ \ \ \  \ \
\int\limits_\R\!\!\pk{
\exists \ttt \in [-S,0)_\delta : Z_u(\ttt)-c_2\ttt>x \text{ or }
\exists t \in [0,S]_\delta : B^*(t)-c_1t>x+(c_1-c_2)\delta
}e^{\frac{\eta x}{t_u}-\frac{(x-c_2\theta_u)^2}{2ut_u}}dx
\\&\leq&
\!\!e^{\frac{\theta_u\eta(\eta-2t^*c_2)}{2(t^*)^2}}
\!\!\!\int\limits_\R\!\!\pk{
\exists \ttt \in [-S,0)_\delta : Z_u(\ttt)-c_2\ttt>x \text{ or }
\exists t \in [0,S]_\delta : B^*(t)-c_1t>x+(c_1-c_2)\theta_u
}e^{\frac{\eta x}{t_u}-\frac{(x-c_2\theta_u)^2}{2ut_u}}dx
\notag\\&\leq&
\!\!e^{\frac{\delta\eta(\eta-2t^*c_2)}{2(t^*)^2}}
\!\!\!\int\limits_\R\!\!\pk{
\exists \ttt \in [-S,0)_\delta : Z_u(\ttt)-c_2\ttt>x \text{ or }
\exists t \in [0,S]_\delta : B^*(t)-c_1t>x
}e^{\frac{\eta x}{t_u}-\frac{(x-c_2\theta_u)^2}{2ut_u}}dx \notag
.}

We estimate the integral in the lower bound. Assume, that BM is defined on
$\R$ (centered Gaussian process with
$\cov(B(t),B(s)) = \frac{|t|+|s|-|s-t|}{2}$).
When $u \rightarrow \IF$ covariance and expectation of
$Z_u(t)-\frac{t\eta}{t^*}$ converge to those of BM,
hence $Z_u(t)-\frac{t\eta}{t^*}$ converges to $B(t)$ for $t<0$
in the sense of convergence of finite-dimensional distributions. Thus
(proof is given in the Appendix)
\bqn{ \label{integral}
& \ &
\int\limits_\R\pk{
\exists \ttt \in [-S,0)_\delta : Z_u(\ttt)-c_2\ttt>x \text{ or }
\exists t \in [0,S]_\delta : B^*(t)-c_1t>x+(c_1-c_2)\delta
}e^{\frac{\eta x}{t_u}-\frac{(x-c_2\theta_u)^2}{2ut_u}}dx
\\&\sim&
\int\limits_\R\pk{
\exists \ttt \in [-S,0)_\delta : B(\ttt)+\zeta\ttt>x \text{ or }
\exists t \in [0,S]_\delta : B(t)-c_1t>x+(c_1-c_2)\delta
}e^{\frac{\eta x}{t^*}}dx, \ \ u \to \IF,
\notag}
with $\zeta = \frac{\eta}{t^*}-c_2$. 
Denote $\xi = \frac{\eta^2}{(t^*)^2}$ and $\kappa = \frac{\eta
\delta(c_1-c_2)}{t^*}$.
We have as $u \to \IF$
\bqny{\label{borell}
& \ &
\int\limits_\R\pk{
\exists \ttt \in [-S,0)_\delta : B(\ttt)+\zeta\ttt>x \text{ or }
\exists t \in [0,S]_\delta : B(t)-c_1t>x+\delta(c_1-c_2)
}e^{\frac{\eta x}{t^*}}dx
\\&=&
\frac{t^*}{\eta}\int\limits_\R\pk{
\exists \ttt \in [-S,0)_\delta :
 B(\ttt)\frac{\eta}{t^*}+\frac{\zeta\eta}{t^*}\ttt>x \text{ or }
\exists t \in [0,S]_\delta : B(t)\frac{\eta}{t^*}-\frac{\eta c_1t}{t^*}>x
+\kappa}e^xdx
\notag\\ &=&
\frac{t^*}{\eta}\int\limits_\R\pk{
\exists \ttt\xi \in [-S\xi,0)_{\delta\xi}: B(\ttt\xi)+\frac{\zeta t^*}{\eta}\ttt\xi>x\text{ or }
\exists t\xi \in [0,S\xi]_{\delta\xi} : B(t\xi)-\frac{\xi c_1t^*t}{\eta}>x
+\kappa }e^xdx
\notag\\ &=&
\frac{t^*}{\eta}\int\limits_\R\pk{
\exists \ttt \in [-S\xi,0)_{\delta\xi}: B(\ttt)+\frac{\zeta t^*}{\eta}\ttt>x \text{ or }
\exists t \in [0,S\xi]_{\delta\xi} : B(t)-\frac{ c_1t^*}{\eta}t>x+\kappa
} e^xdx
\notag\\  &=:&
\frac{t^*}{\eta} I(S)
\notag.}
Similarly to the proof of $\eqref{finitness_J(S)}$ we have
\bqn{\label{finitness_I(S)}
\limit{S}I(S) \in (0,\IF).
}

Denote
$$\hat{d}(t) = \mathbb{I}(t<0)\frac{\zeta t^*t}{\eta} -
\mathbb{I}(t\geq 0)(\frac{c_1t^*t}{\eta}+\kappa).$$
We have
\bqny{ I(S) =
\int\limits_\R\pk{\exists t \in [-S\xi,S\xi]_{\xi\delta}
:B(t)+\hat{d}(t)>x
}e^xdx
&=&
\E{\exp\left(\sup\limits_{t \in [-S\xi,S\xi]_{\delta\xi}} (B(t)+\hat{d}(t))\right)} \\ &=&
\E{\sup\limits_{t \in [-\frac{S\xi}{2},\frac{S\xi}{2}]_\frac{\delta\xi}{2}}
e^{\sqrt 2 B(t)-|t|+\hat{d}(2t)+|t|}}.
}
We have that
$$ \hat{d}(2t)+|t|  =
\mathbb{I}(t<0)\frac{(q_2c_1+c_2q_1-2q_2c_2)t}{c_1q_2-q_1c_2}  +
\mathbb{I}(t\geq 0)\Big(\frac{(2c_1q_1-c_1q_2-q_1c_2)t}{c_1q_2-q_1c_2}-
\delta\frac{(c_1q_2-q_1c_2)(c_1-c_2)}{q_2-q_1}\Big)
=d_\delta(t).$$
Since  $\frac{\delta\xi}{2}  = \gamma$
by \eqref{finitness_I(S)} we have
$$\E{\sup\limits_{t \in [-\frac{S\xi}{2},\frac{S\xi}{2}]_\frac{\delta\xi}{2}}
e^{\sqrt 2 B(t)-|t|+\hat{d}(2t)+|t|}}  \rightarrow \mathcal{H}^{d_\delta}
_\gamma \in (0,\IF), \ \ \
S \to \IF. $$

Thus we conclude that as $u \to \IF$ and then $S \to \IF$
\bqn{\label{lowbound}
\pk{\sup\limits_{t \in \Delta}Z(t)>\sqrt u} \geq \mathcal{H}^{d_\delta}
_\gamma \PH(\mathbb{D}_{1/2}\sqrt u)(1+o(1))
.}
For the same reasons we have that as $u \to \IF$ and then $S \to \IF$
\bqn{\label{upbound}
\pk{\sup\limits_{t \in \Delta}Z(t)>\sqrt u} \leq
\mathcal{H}^{d}_\gamma e^{\frac{\delta\eta(\eta-2t^*c_2)}{2(t^*)^2}}
\PH(\mathbb{D}_{1/2}\sqrt u)
(1+o(1))
.}
Finally we can compare the approximations in \eqref{sums}.
From \eqref{lowbound} and \eqref{asympsum2} it follows
that $\sum\limits_{j=1}^{N_u}p_{j,S,u}$ is negligible
as $u \to \IF$ and then $S \to \IF$.
By the same arguments $\sum\limits_{j=-N_u}^{-2}p_{j,S,u}$
is also negligible. Thus the claim follows
from \eqref{lowbound} and \eqref{upbound}.
\\
\\
\emph{Case $H<1/2$.}
First of all we prove the lower bound in \eqref{small}.
We have with $t_u^- = t_u-\delta/u$ and $w = u^{1-H}$
($V_1(t)$ and $V_2(t)$ are defined in \eqref{videf})
$$
\bar\psi_{\delta,H}(u)
\ge
\pk{ \sup\limits_{t\in \{t_u^-,t_u\}}Z_H(t)>w }
=
\pk{V_1(t_u)>w}+\pk{V_2(t_u^-)>w}
-\pk{V_1(t_u)>w, V_2(t_u^-)>w}.
$$
We have by the ratio $\PH(x)\sim \frac{1}{x\sqrt{2\pi}}e^{-x^2/2}, \ x \to
\IF$ and Taylor's theorem as $u \to \IF$
\bqny{
\pk{V_1(t_u)>w}
 =
\PH(u^{1-H}\frac{q_1+c_1t_u}{t_u^H}) \sim
\PH(\md u^{1-H})\exp(-\frac{\theta_u w_1'(t^*)}{2}u^{1-2H})
,}
where $w_1(t)$ is defined in \eqref{ww}.
By the same arguments
\bqny{
\pk{V_2(t_u^-)>w} \sim
\PH(\md u^{1-H})\exp(-\frac{-(\delta-\theta_u) w_2'(t^*)}{2}u^{1-2H}),
\ \ \ \ u \to \IF.
}
By Lemma 2.3 in \cite{PickandsA} we conclude that as $u \to \IF$
$$\pk{V_1(t_u)>w,
V_2(t_u^-)>w}
= o(\pk{V_1(t_u)>w}) $$
and hence as $u \to \IF$
$$\pk{\sup\limits_{t\in \{t_u^-,t_u\}}Z_H(t)>w} =
\PH(\md u^{1-H})\left(\exp(-\frac{\theta_u w_1'(t^*)}{2}u^{1-2H})+
\exp(-\frac{-(\delta-\theta_u) w_2'(t^*)}{2}u^{1-2H})
\right)(1+o(1)).$$
Notice that $w_1'(t^*)>0$ and $w_2'(t^*)<0$.
We have as $u \to \IF$ (recall, $B$ is defined in \eqref{ww})
\bqn{\label{remark1}
\exp(-\frac{\theta_u w_1'(t^*)}{2}u^{1-2H}) +
\exp(-\frac{-(\delta-\theta_u) w_2'(t^*)}{2}u^{1-2H})
\geq
2e^{-Bu^{1-2H}}(1+o(1))
.}
Hence the lower bound in \eqref{small} is established.\\ 

Next we prove the upper bound. We have as $u \to \IF$
\bqn{\label{sm}
\bar{\psi}_{\delta,H}(u) &\leq&
\pk{\exists t \in G(\delta/u),t<t^*:V_2(t)>w}+
\pk{\exists t \in G(\delta/u),t\geq t^*:V_1(t)>w}
\notag\\&\sim&
\pk{V_2(t_u^-)>w}+
\pk{V_1(t_u)>w},
}
proof of the last line above is in the Appendix.
Similarly to the proof of the lower bound we have
$$\pk{V_2(t_u^-)>w}
+ \pk{V_1(t_u)>w}\leq \pk{V_1(t^*)>w}(1+o(1))
\sim \PH(\md u^{1-H}), \ \ u \to \IF$$
establishing the claim.\\

\textbf{Case (3).}
Assume that $i=1$, case $i=2$ follows by the same arguments. We have
\bqn{\label{ineqH4}
\pk{\sup\limits_{t \in [t^*,\IF)_\frac{\delta}{u}} V_1(t)>u^{1-H}}
&\leq& \bar{\psi}_{\delta,H}(u)
\notag\\ &\leq&
\pk{\sup\limits_{t \in [t^*,\IF)_\frac{\delta}{u}} V_1(t)>u^{1-H}}
+ \pk{\sup\limits_{t \in [0,t^*)_\frac{\delta}{u}} V_2(t)>u^{1-H}}
.}

\emph{Case $H \geq 1/2$.}
It follows from Theorem 1 in \cite{Piterbargdiscrete} that
$$\pk{\sup\limits_{t \in [t^*,\IF)_\frac{\delta}{u}} V_1(t)>u^{1-H}}
\sim \frac{1}{2}\psi^{(1)}_{\delta,H}(u), \ \ \ \ u \to \IF,$$
where $\psi^{(1)}_{\delta,H}(u)$ is defined in \eqref{psi1def}.
Since $H \ge \frac{1}{2}$ in view of
Remark \ref{remcompasymp} for some positive constant $C$
that does not depend on $u$
$$\psi^{(1)}_{\delta,H}(u) \ge  C\psi^{(1)}_{0,H}(u), \ \ \ u \to \IF.$$
It follows from the proof of Theorem 3.1, case (4), $H \geq 1/2$ in \cite{Lanpeng2BM} that
$$\pk{\sup\limits_{t \in [0,t^*)} V_2(t)>u^{1-H}}
= o(\psi^{(1)}_{0,H}(u)), \ \ \ u \to \IF.$$
Hence \eqref{ineqH4} yields that
$$\bar{\psi}_{\delta,H}(u) \sim \frac{1}{2}\psi^{(1)}_{\delta,H}(u), \ \ \ \ u \to \IF$$
establishing the claim by \eqref{one-dimensional_theorem}.\\

\emph{Case $H<1/2$.}
From proof of Theorem \ref{H1}, case $H<1/2$ follows, that
$$\pk{\sup\limits_{t \in [t^*,\IF)_\frac{\delta}{u}} V_1(t)>u^{1-H}} \sim
\frac{1}{2}\psi^{(1)}_{\delta,H}(u), \ \ u \to \IF.$$
Notice that
\bqny{
\pk{\sup\limits_{t \in [0,t^*)_\frac{\delta}{u}} V_2(t)>u^{1-H}}
&\leq&
\pk{V_2(t_u^-)>u^{1-H}}+\frac{ut^*}{\delta}
\sup\limits_{t \in [0,t^*-\delta/u)_\frac{\delta}{u}}
 \pk{V_2(t)>u^{1-H}}
\\ &=&
\pk{V_2(t_u^-)>u^{1-H}} +
\frac{ut^*}{\delta} \pk{V_2(t_u^--\delta/u)>u^{1-H}}
\\ &=&
\PH(\frac{u^{1-H}(q_2+c_2t_u^-)}{t_u^H}) +
\frac{ut^*}{\delta}\PH(\frac{u^{1-H}(q_2+c_2(t_u^--\delta/u))}
{(t^-_u-\delta/u)^H})
.}
We have that
\bqny{
\PH(\frac{u^{1-H}(q_2+c_2(t_u^--\delta/u))}
{(t_u^-\delta/u)^H})
= \PH(\md u^{1-H})\exp(u^{1-2H}\frac{(2\delta-\theta_u)w_2'(t^*)}{2})(1+o(1)), \ \ u \to \IF.
}
Since $H<1/2$ and $w_2'(t^*)<0$ it follows from
\eqref{one-dimensional_theorem} that the expression above
equals $o(\psi^{(1)}_{\delta,H}(u))$ as $u \to \IF$.
Hence from \eqref{ineqH4} it follows
$$\bar{\psi}_{\delta,H}(u) \sim \frac{1}{2}\psi^{(1)}_{\delta,H}(u), \ \ u \to \IF$$
and the claim is established.
\QED
\\
\\
\textbf{Proof of Remark \ref{rem2}.}
Consider a sequence $\{u_n\}_{n\in \mathbb{N}}$ such $u_n \to \IF$
and for all $n$ $t^*\in G(\delta/u_n)$.
From the proof of Theorem \ref{maintheo} case (2), $H<1/2$ it follows, that
$$\bar{\psi}_{\delta,H}(u_n) =
\PH(\md u_n^{1-H})(1+o(1)), \ \ \ n \to \IF.$$
We choose a sequence
$\{v_n\}_{n\in \mathbb{N}}$ such $v_n \to \IF$ and for all $n$
$t^*-\frac{\delta w'_2(t^*)}{v_n(w'_1(t^*)-w'_2(t^*))} \in G(\delta/v_n)$.
For such sequence inequality in \eqref{remark1} becomes
equality, hence
\bqny{
\bar{\psi}_{\delta,H}(v_n) \sim
2e^{-Bv_n^{1-H}}\PH(\md v_n^{1-H})
, \ \ \ \ \ n \to \IF.
\ \ \ \ \ \ \ \ \ \ \ \ \ \ \hfill \Box}

\textbf{Proof of Theorem \ref{H1}.}
When $H=\frac{1}{2}$ the assertion of the theorem follows from
the results in \cite{HashorvaGrigori} and \cite{Piterbargdiscrete}.
\\
\\
\emph{Case $H > 1/2$.}
For large $u>0$
$$
\pk{\exists t \ge 0 : \inf\limits_{s\in [t,t+u^{\frac{2H-1}{2H}}]}
(B_H(s)-cs)>u}
\le
\pk{\sup\limits_{t\in G(\delta)}(B_H(t)-ct)>u}
\le
\pk{\sup\limits_{t\ge 0}(B_H(t)-ct)>u}
.$$
In a view of Remark 3.4 in \cite{HashorvaJiDebickiParisianSelfSimilar}
$$
\pk{\exists t \ge 0 : \inf\limits_{s\in [t,t+u^{\frac{2H-1}{2H}}]}
(B_H(s)-cs)>u}
\sim
\pk{\sup\limits_{t\ge 0}(B_H(t)-ct)>u}, \ \ \ u \to \IF
$$
implying
$$\pk{\sup\limits_{t\in G(\delta)}(B_H(t)-ct)>u}
\sim
\pk{\sup\limits_{t\ge 0}(B_H(t)-ct)>u}, \ \ \ u \to \IF.$$
The asymptotic of the last probability above is given e.g. in
Proposition 2.1 in \cite{Lanpeng2BM}, thus the claim follows.
\\
\\
\emph{Case $H<1/2$.} By the self-similarity of fBM we have
\bqn{\label{psidef}
\psi_{\delta,H}(u) :=
\pk{\exists t \in G(\delta): B_H(t)>u+ct} &=&
\pk{\exists t \in G(\frac{\delta}{u}):\frac{B_H(t)}{1+ct}>u^{1-H}}
\\&=:& \pk{\exists t \in G(\frac{\delta}{u}):V(t)>u^{1-H}}
\notag.}
Note, that the variance of $V(t)$ achieves its unique maxima at
$t_0=\frac{H}{c(1-H)}$. 
As shown in the Appendix
\bqn{\label{negl}
\psi_{\delta,H}(u) \sim \pk{\sup\limits_{t \in I(t_0)}V(t)>u^{1-H}}
:= \varsigma(u), \ \ \ \ u \to \IF,
}
with $I(t_0) = (-1/\sqrt u+t_0,1/\sqrt u+t_0)$.
Next (proof see in the Appendix)
\bqn{\label{bonf}
\varsigma(u) \sim \sum\limits_{t \in I(t_0)}\pk{V(t)>u^{1-H}}, \  \ u \to \IF.
}

We approximate the sum above. We have as $u \to \IF$
with $\hu = \frac{u^{1-H}c^H}{H^H(1-H)^{1-H}}$
\bqny{
\sum\limits_{t \in I(t_0)}\pk{V(t)>u^{1-H}} =
\sum\limits_{t \in I(t_0)}\PH(u^{1-H}\frac{1+ct}{t^H}) \sim
\sum\limits_{t \in I(t_0)}\frac{1}{\sqrt {2\pi} \hu}e^{-\frac{1}{2}(u^{1-H}\frac{1+ct}{t^H})^2}
.}
Define
\bqn{\label{ff}
f_H(t) = \frac{(1+ct)^2}{t^{2H}}.
}
We have that $f_H'(t_0) = 0$ and
$f_H''(t_0)=\frac{2c^{2+2H}(1-H)^{2H+1}}{H^{2H+1}}>0$.
We write as $u \to \IF$
\bqn{\label{comput}
\sum\limits_{t \in I(t_0)}\frac{1}{\sqrt {2\pi} \hu}e^{-\frac{1}{2}(u^{1-H}\frac{1+ct}{t^H})^2}
&=&
\frac{1}{\sqrt {2\pi} \hu}e^{-\hu^2/2}\sum\limits_{t \in I(t_0)} e^{-\frac{1}{2}u^{2-2H}\left(\frac{(1+ct)^2}{t^{2H}} -
\frac{(1+ct_0)^2}{t_0^{2H}}\right)}
\notag \\&\sim&
\PH(\hu)\sum\limits_{t \in I(t_0)}e^{-\frac{1}{2}u^{2-2H}\frac{f_H''(t_0)}{2}(t-t_0)^2}.
}

The proof of the last line above is given in the Appendix.
Next (set $F = \frac{f_H''(t_0)}{4} = \frac{c^{2+2H}(1-H)^{2H+1}}{2H^{2H+1}}$)
\bqny{\sum\limits_{t \in I(t_0)}e^{-\frac{1}{2}u^{2-2H}\frac{f_H''(t_0)}{2}(t-t_0)^2}
&\sim&
2\sum\limits_{t \in (0,u^{-1/2})_{\delta/u}}e^{-Fu^{2-2H}t^2}
\\ &=&
2\sum\limits_{tu^{1-H} \in (0,u^{1/2-H})_{\delta u^{-H}}}e^{-F(tu^{1-H})^2}
\\ &=&
\frac{2u^H}{\delta} \big(\delta u^{-H}
\sum\limits_{t \in (0,u^{1/2-H})_{\delta u^{-H}}}e^{-Ft^2}\big)
\\ &\sim&
\frac{2u^H}{\delta\sqrt F}\int\limits_0^\IF e^{-Ft^2}d(\sqrt{F}t)
\\ &=&
\frac{\sqrt{\pi}u^H}{\delta\sqrt F}, \ \ u \to \IF.}
Combining the line above with \eqref{comput} we have
\bqn{\label{fincomput}
\varsigma(u) \sim \PH(\frac{u^{1-H}c^H}{H^H(1-H)^{1-H}})\frac{\sqrt{2\pi} H^{H+1/2} u^H}{\delta c^{H+1}(1-H)^{H+1/2}}, \ \ \ \ u \to \IF
,}
and by \eqref{negl} and \eqref{bonf}
the claim is established. \QED
\\
\\
\textbf{Proof of Remark \ref{prop}.} Assume, that $(q_1,c_1) \ge (q_2,c_2)$ in the alphabetical order.
Then for large $u$  $q_1u+c_1t \geq q_2u+c_2t$ for all $t \in [0,T]$, hence for large $u$
\bqny{
\bar{\zeta}_H(u) &=& \pk{\exists t \in [0,T]: B_H(t)-c_1t>q_1u,B_H(t)-c_2t>q_2u}\\
&=& \pk{\exists t \in [0,T]: B_H(t)>\max(c_1t+q_1u,c_2t+q_2u)}\\
&=& \pk{\exists t \in [0,T]: B_H(t)>c_1t+q_1u}
.}
If
$(q_2,c_2) > (q_1,c_1)$ in the alphabetical order, then by the same arguments
\bqny{
\bar{\zeta}_H(u) = \pk{\exists t \in [0,T]: B_H(t)>c_2t+q_2u}
.}
Consequently, $\bar{\zeta}_H(u)$ for large $u$ always coincides with one of the single-dimensional probabilities, namely, with the smallest one.
For $H=1/2$ the claim follows by \cite{DeM15}.
For $H\neq 1/2$ Theorem 2.1 in \cite{SumFBMDebicki}
establishes the proof. \QED

\section{Appendix}

\textbf{Proof of \eqref{borH}.}
We shall prove that $V_1(t)$ is a.s. bounded on $[0,\IF)$. By Chapter 4, p. 31 in
\cite{20lectures} it is equivalent to
$$\pk{V_1(t) \text{ is bounded for } t\ge 0}>0.$$
We have as $u \to \IF$
\bqny{
\pk{\sup\limits_{t\ge 0}V_1(t)\le u} &=& 1-\pk{\sup\limits_{t\ge 0}V_1(t)> u}
\to 1
}
by Proposition 2.1 in \cite{Lanpeng2BM}.
Thus $V_1(t)$ is bounded a.s.
\\
\\
Notice, that the variance $v(t)$ of
$V_1(t)$ achieves its unique maxima at $t_1$.
Denote
$$m = \max\limits_{t \in [0,t_1-\ve]\cap[t_1+\ve,\IF)}v(t), \ \ \ M = \E{\sup\limits_{t \in [0,t_1-\ve]\cap[t_1+\ve,\IF)}V_1(t)}.$$
By Borell-TIS inequality (see Lemma 5.3 in \cite{Lanpeng2BM}) we
have that
$M < \IF$ and for all $u$ large enough we have
\bqny{
\pk{\exists t \notin [t_1-\ve,t_1+\ve]: V_1(t)>u^{1-H}}
 \leq
e^{-\frac{(u^{1-H}-M)^2}{2m}}
.}
From Theorem \ref{H1} and inequality $m<v(t_1)$ it follows, that
$$e^{-\frac{(u^{1-H}-M)^2}{2m}} = o(\psi^{(1)}_{\delta,H}(u))
, \ \ \ u \to \IF$$
and thus \eqref{borH} holds.
\QED
\\
\\
\textbf{Proof of \eqref{informinterv}.}
We have
\bqn{
\pk{\exists t \in [t^*-\frac{\ln u}{\sqrt u}, t^*+\frac{\ln u}{\sqrt u}]_{\frac{\delta}{u}} : Z(t)> \sqrt u} \leq \bp_\delta(u) \leq
\pk{\exists t \in [t^*-\frac{\ln u}{\sqrt u}, t^*+\frac{\ln u}{\sqrt u}]_{\frac{\delta}{u}} : Z(t)> \sqrt u}
\notag\\ \label{00}
+ \pk{\exists t \in [0, t^*-\frac{\ln u}{\sqrt u}]_{\frac{\delta}{u}} : Z(t)> \sqrt u}
+ \pk{\exists t \in [t^*+\frac{\ln u}{\sqrt u}, \IF)_{\frac{\delta}{u}} : Z(t)> \sqrt u}
.}
We have for some fixed $\ve>0$
\bqny{
\pk{\exists t \in [t^*+\frac{\ln u}{\sqrt u}, \IF)_{\frac{\delta}{u}} : Z(t)> \sqrt u}
&\leq&
\pk{\exists t \in [t^*+\frac{\ln u}{\sqrt u}, t^*+\ve)_{\frac{\delta}{u}} : Z(t)> \sqrt u}
\\&+&
\pk{\exists t \in [t^*+\ve, \IF)_{\frac{\delta}{u}} : Z(t)> \sqrt u}.
}
Thus by Borell-TIS inequality ($Z(t)$ is bounded a.s. by the same
arguments as in the proof of \eqref{borH}),
\bqny{
\pk{\exists t \in [t^*+\ve, \IF)_{\frac{\delta}{u}} : Z(t)> \sqrt u} =
o(\pk{Z(t_u)}>\sqrt u), \ \ \ \ u \to \IF.
}
Also
\bqny{
\pk{\exists t \in [t^*+\frac{\ln u}{\sqrt u}, t^*+\ve)_{\frac{\delta}{u}} : Z(t)> \sqrt u}
&\leq&
\frac{2\ve u}{\delta}\sup\limits_{t \in [t^*+\frac{\ln u}{\sqrt u}, t^*+\ve)
_{\frac{\delta}{u}}}\pk{Z(t)>\sqrt u}
\\&\leq&
\frac{2\ve u}{\delta}\pk{Z(t^*+\frac{\ln u}{\sqrt u})>\sqrt u}
\notag\\ &=&
\frac{2\ve u}{\delta}\PH(\frac{\sqrt u(c_1(t^*+\frac{\ln u}{\sqrt u})+q_1)}
{\sqrt{t^*+\frac{\ln u}{\sqrt u}}})
\notag\\&=&
o(\pk{Z(t_u)}>\sqrt u), \ \ \ u \to \IF
.\notag}
Thus we conclude, that
\bqny{
\pk{\exists t \in [t^*+\frac{\ln u}{\sqrt u}, \IF)
_{\frac{\delta}{u}} : Z(t)> \sqrt u} = o(\pk{Z(t_u)>\sqrt u}), \ \ u \to \IF.
}
By the same arguments
\bqny{
\pk{\exists t \in [0,t^*-\frac{\ln u}{\sqrt u})
_{\frac{\delta}{u}} : Z(t)> \sqrt u} = o(\pk{Z(t_u)>\sqrt u}), \ \ u \to \IF.
}
Hence \eqref{informinterv} follows from
\eqref{00} and two expressions above. \QED
\\
\\
\textbf{Proof of \eqref{derivate}.} Recall, that $a(t) = \frac{(1+\mu t)^2}{2t}$.
We have
\bqny{& \ &
|\sum\limits_{j=1}^{N_u} e^{-v^2(\frac{(1+\mu c_{j,S,v})^2}{2c_{j,S,v}}-\frac{(1+\mu t^*)^2}{2t^*})} -  \sum\limits_{j=1}^{N_u} e^{-v^2a'(t^*)(\frac{\theta_u+jS}{u})}|
\\&=&
|\sum\limits_{j=1}^{N_u}(e^{-v^2(a(c_{j,S,v})-a(t^*))} -
e^{-v^2a'(t^*)(\frac{\theta_u+jS}{u})})| \\
&=&
|\sum\limits_{j=1}^{N_u}e^{-v^2a'(t^*)(\frac{\theta_u+jS}{u})}(
e^{-v^2\big(a(c_{j,S,v})-a(t^*)-a'(t^*)(c_{j,S,v}-t^*)\big)}-1)| \\
&\leq&
u^{1/10}\sup\limits_{1\leq j\leq u^{1/10}}e^{-v^2a'(t^*)(\frac{\theta_u+jS}{u})}
\sup\limits_{1\leq j\leq u^{1/10}}|
e^{-v^2\big(a(c_{j,S,v})-a(t^*)-a'(t^*)(c_{j,S,v}-t^*)\big)}-1|
+CN_ue^{-u^{1/10}C_1}
\\&\leq&
Cu^{1/10}\sup\limits_{1\leq j\leq u^{\frac{1}{10}}}|e^{-v^2\big(a(c_{j,S,v})-a(t^*)-a'(t^*)(c_{j,S,v}-t^*)
\big)}-1|
+CN_ue^{-u^{1/10}C_1}.}

Notice that $a''(t^*) = \frac{1}{(t^*)^3}>0$.
When $u$ is large enough and $j \leq u^{\frac{1}{10}}$
\bqny{
|e^{-v^2\big(a(c_{j,S,v})-a(t^*)-a'(t^*)(c_{j,S,v}-t^*)\big)}-1| \leq
|e^{-v^2a''(t^*)(c_{j,S,v}-t^*)^2}-1| =
|e^{-Cv^2(\frac{\theta_u+jS}{u})^2}-1| \leq u^{-7/10}
.}
Thus
\bqny{
Cu^{1/10}\sup\limits_{1\leq j\leq u^{\frac{1}{10}}}|e^{-v^2\big(a(c_{j,S,v})-a(t^*)-a'(t^*)(c_{j,S,v}-t^*)\big)}-1|
< u^{-1/2}.}
Hence as $u \to \IF$
\bqn{\label{to0}
|\sum\limits_{j=1}^{N_u} e^{-v^2(\frac{(1+\mu c_{j,S,v})^2}{2c_j}-\frac{(1+\mu t^*)^2}{2t^*})} -  \sum\limits_{j=1}^{N_u} e^{-v^2a'(t^*)(\frac{\theta_u+jS}{u})}|<u^{-1/2}+CN_ue^{-u^{1/10}C_1}
\to 0.}
As was shown in the proof of Theorem \eqref{maintheo}, case (2)
\bqny{
\sum\limits_{j=1}^{N_u} e^{-v^2a'(t^*)(\frac{\theta_u+jS}{u})} \sim
e^{-\omega(\theta_u+S)}
\ge
e^{-\omega(\delta+S)} >0,  \ \ \ \ \ u \to \IF .
}
Thus the claim follows by the line above and \eqref{to0}.
\QED
\\
\\

\textbf{Proof of \eqref{gausint}.}
We find the distribution law of $X_u(\ttt) := \{B(\ttt)|B(ut_u) = \eta u-x\}, \ttt \in [ut_u-S,ut_u]$.
Let $ut_u-S \leq \sss \leq \ttt \leq ut_u$ and  $N_1, N_2$
be independent standard Gaussian random variables
independent of $B(ut_u)$. Then
$$(B(\sss),B(\ttt),B(ut_u)) = (\sqrt{\sss-\frac{\sss^2}{\ttt}}N_1+\frac{\sss\sqrt{ut_u-\ttt}}{\sqrt{\ttt ut_u}}N_2+\frac{\sss}{ut_u}B(ut_u),
\sqrt{\ttt-\frac{\ttt^2}{ut_u}}N_2+\frac{\ttt}{ut_u}B(ut_u),B(ut_u)).$$
By the formula above, we have
$$\E{X_u(\ttt)} = \frac{\ttt(\eta u-x)}{ut_u},
 \ \ \ \ \ \ \cov(X_u(\sss),X_u(\ttt)) = \sss-\frac{\sss\ttt}{ut_u}.$$

For the process $$Y_u(\ttt)+\eta u -x  = X_u(\ttt), \ \ \ \ttt \in[-S+ut_u,ut_u]$$
we have
$$\E{Y_u(\ttt)} = \frac{(\eta u-x)(\ttt-ut_u)}{ut_u}, \quad
\cov (Y_u(\sss),Y_u(\ttt)) = \sss-\frac{\sss\ttt}{ut_u}.$$
Thus
\bqny{
& \ &
\PP\{
\exists \ttt \in [ut_u-S,ut_u)_\delta : B(\ttt)>q_2u+c_2ut_u+c_2(\ttt-ut_u)
\\ & \ & \text{ or }
\exists t \in [ut_u,ut_u+S]_\delta : B^*(t-ut_u)+\eta u -x>q_1u+c_1ut_u+c_1(t-ut_u)
 |B(ut_u) = \eta u-x\}
\\&=&
\pk{
\exists \ttt \in [ut_u-S,ut_u)_\delta : Y_u(\ttt)>x+c_2\theta_u+c_2(\ttt-ut_u)
\text{ or }
\exists t \in [0,S]_\delta : B^*(t)>x+c_1\theta_u+c_1t
}
\\&=&
\pk{
\exists \ttt \in [-S,0)_\delta : Z_u(\ttt)-c_2\ttt>x+c_2\theta_u \text{ or }
\exists t \in [0,S]_\delta : B^*(t)-c_1t>x+c_1\theta_u
}.
}
Recall, that $Z_u(\ttt)$ is a Gaussian process with expectation
and covariance defined below ($\sss \leq \ttt$):
$$\E{Z_u(\ttt)} = \frac{\eta u-x}{ut_u}\ttt,
 \ \ \ \ \ \ \cov (Z_u(\sss),Z_u(\ttt)) = \frac{-\sss\ttt}{ut_u}-\ttt.
$$
\QED

\textbf{Proof of \eqref{integral}.}
First of all we show that with $\bar \delta = (c_1-c_2)\delta$
\bqn{\begin{split}\label{apborel}
\int\limits_\R\pk{
\exists \ttt \in [-S,0)_\delta : Z_u(\ttt)-c_2\ttt>x \text{ or }
\exists t \in [0,S]_\delta : B(t)-c_1t>x+\bar\delta
}e^{\frac{\eta x}{t_u}-\frac{(x-c_2\theta_u)^2}{2ut_u}}dx\\
=\int\limits_{-M}^M\pk{
\exists \ttt \in [-S,0)_\delta : Z_u(\ttt)-c_2\ttt>x \text{ or }
\exists t \in [0,S]_\delta : B(t)-c_1t>x+\bar\delta
}e^{\frac{\eta x}{t^*}}dx+B_{M,v},
\end{split}}
where $B_{M,v} \to \IF$ when $u\to \IF$ and then $M\to \IF$. We have
\bqny{& \ &
\Big|\int\limits_\R\pk{
\exists \ttt \in [-S,0)_\delta : Z_u(\ttt)-c_2\ttt>x \text{ or }
\exists t \in [0,S]_\delta : B(t)-c_1t>x+\bar\delta
}e^{\frac{\eta x}{t_u}-\frac{(x-c_2\theta_u)^2}{2ut_u}}dx
\\& \ & - \ \
\int\limits_{-M}^{M}\pk{
\exists \ttt \in [-S,0)_\delta : Z_u(\ttt)-c_2\ttt>x \text{ or }
\exists t \in [0,S]_\delta : B(t)-c_1t>x+\bar\delta
}e^{\frac{\eta x}{t^*}}dx\Big|
\\&\leq&
\Big|\int\limits_{-M}^{M}\pk{
\exists \ttt \in [-S,0)_\delta : Z_u(\ttt)-c_2\ttt>x \text{ or }
\exists t \in [0,S]_\delta : B(t)-c_1t>x+\bar\delta
}(e^{\frac{\eta x}{t_u}-\frac{(x-c_2\theta_u)^2}{2ut_u}}
-e^{\frac{\eta x}{t^*}})dx
\Big|
\\& \ & + \ \
\int\limits_{|x|>M}\pk{
\exists \ttt \in [-S,0)_\delta : Z_u(\ttt)-c_2\ttt>x \text{ or }
\exists t \in [0,S]_\delta : B(t)-c_1t>x+\bar\delta
}e^{\frac{\eta x}{t_u}-\frac{(x-c_2\theta_u)^2}{2ut_u}}dx
\\&=:&
|I_1|+I_2
.}

By Borell-TIS inequality for large $u$ and $x>0$
\bqny{ & \ &
\pk{\exists \ttt \in [-S,0)_\delta : Z_u(\ttt)-c_2\ttt>x \text{ or }
\exists t \in [0,S]_\delta : B(t)-c_1t>x+\bar\delta}
\\&\le&
\pk{\exists t \in [-S,0] : Z_u(t)-\E{Z_u(t)}>x}
+
\pk{\exists t \in [0,S] : B(t)>x}
\\&\le&
e^{-x^2/C}.
}
Thus as $u \to \IF$
\bqny{
I_2 \leq \int\limits_{x>M}
e^{-\frac{x^2}{C}+\frac{\eta x}{t_u}-\frac{(x-c_2\theta_u)^2}{2ut_u}}
dx
+\int\limits_{-\IF}^{-M}e^{\frac{\eta x}{2t^*}}dx
 \to 0, \ \ \ M  \to \IF.
}
For $I_1$ we have for large $u$
\bqny{
|I_1| \leq
\int\limits_{-M}^{M}
e^{-\frac{x^2}{C}+\frac{\eta x}{t^*}}(e^{-
\frac{x\eta\theta_u }{ut^*t_u}-\frac{(x-c_2\theta_u)^2}{2ut_u}}-1)dx
.}
Let $u\geq M^3$. For such $u$
$$\sup\limits_{x \in [-M,M]}|e^{-
\frac{x\eta\theta_u }{ut^*t_u}-\frac{(x-c_2\theta_u)^2}{2ut_u}}-1| \leq C\frac{1}{M},$$
hence for $u\geq M^3$
\bqny{& \ &
\int\limits_{-M}^{M}
e^{-\frac{x^2}{C}+\frac{\eta x}{t^*}}(e^{-
\frac{x\eta\theta_u }{ut^*t_u}-\frac{(x-c_2\theta_u)^2}{2ut_u}}-1)dx
\leq
\frac{C}{M}\int\limits_\R
e^{-\frac{x^2}{C}+\frac{\eta x}{t^*}}
=
\frac{C_1}{M}
.}
Thus we proved that
$$\limit{M}\limit{u}(|I_1|+I_2)=0$$
and \eqref{apborel} holds.
Since for $t\in [-S,0] \ Z_u(t)-\frac{\eta}{t^*}t$ converges to $B(t)$
in the sense of convergence of finite-dimensional distributions we have
(recall, $\zeta = \frac{\eta}{t^*}-c_2$)
\bqny{ & \ &
\int\limits_{-M}^M\pk{
\exists \ttt \in [-S,0)_\delta : Z_u(\ttt)-c_2\ttt>x \text{ or }
\exists t \in [0,S]_\delta : B(t)-c_1t>x+\bar\delta
}e^{\frac{\eta x}{t^*}}dx
\\ &\to&
\int\limits_{-M}^M\pk{
\exists \ttt \in [-S,0)_\delta : B(t)+\zeta\ttt>x \text{ or }
\exists t \in [0,S]_\delta : B(t)-c_1t>x+\bar\delta
}e^{\frac{\eta x}{t^*}}dx, \ \ \ \ \ u \to \IF
.}
By the monotone convergence theorem
\bqny{\limit{M}\int\limits_{-M}^M\pk{
\exists \ttt \in [-S,0)_\delta : B(\ttt)+\zeta\ttt>x \text{ or }
\exists t \in [0,S]_\delta : B(t)-c_1t>x+\bar\delta
}e^{\frac{\eta x}{t^*}}dx\\
=
\int\limits_\R\pk{
\exists \ttt \in [-S,0)_\delta : B(\ttt)+\zeta\ttt>x \text{ or }
\exists t \in [0,S]_\delta : B(t)-c_1t>x+\bar\delta
}e^{\frac{\eta x}{t^*}}dx,}
and the claim is established.
\QED
\\
\\

\textbf{Proof of \eqref{sm}.} First we prove
\bqn{\label{app}
\pk{\exists t \in G(\delta/u),t\geq t_u+\delta/u:V_1(t)>u^{1-H}}
= o\left(\pk{V_1(t_u)>u^{1-H}}\right)
, \ \ u \to \IF.
}
Fix some $\ve>0$. By Borell-TIS inequality
 as $u \to \IF$
\bqn{\label{o1}
\pk{\exists t \in G(\delta/u),t\geq t^*+\ve:V_1(t)>u^{1-H}}
=
o\left(\PH(u^{1-H}\frac{c_1t_u+q_1}{t_u^H})\right).
}
We have as $u \to \IF$
\bqny{
\pk{\exists t \in G(\delta/u), t_u+\delta/u\leq t \leq t^*+\ve :
V_1(t)>u^{1-H}}
&\leq&
\frac{\ve u}{\delta}
\sup\limits_{t \in G(\delta/u), t_u+\delta/u\leq t \leq t^*+\ve}
\pk{V_1(t)>u^{1-H}}
\\ &\leq&
\frac{\ve u}{\delta}\pk{V_1(t_u+\delta/u)>u^{1-H}}
\\ &=&
\frac{\ve u}{\delta}\PH(u^{1-H}\frac{c_1(t_u+\delta/u)+q_1}
{(t_u+\delta/u)^H})
\\&\sim&
\frac{\ve u}{\delta}
\PH(u^{1-H}\frac{c_1t_u+q_1}{t_u^H})
\exp(-\frac{w'_1(t^*)\delta}{2}u^{1-2H})
\\ &=&
o\left(\PH(u^{1-H}\frac{c_1t_u+q_1}{t_u^H})\right)
.}
Hence
\bqny{
\pk{\exists t \in G(\delta/u), t_u+\delta/u\leq t \leq t^*+\ve :V_1(t)
>u^{1-H}}
=
o\left(\PH(u^{1-H}\frac{c_1t_u+q_1}{t_u^H})\right), \ \ u \to \IF.
}
By \eqref{o1} with the the line above we establish \eqref{app}.
By the same arguments
$$
\pk{\exists t \in G(\delta/u),t< t_u^-:V_2(t)>u^{1-H}}
= o\left(\pk{V_2(t_u^-)>u^{1-H}}\right)
, \ \ u \to \IF.
$$
Combination of \eqref{app} and the line above establishes the claim.\QED
\\
\\
\textbf{Proof of \eqref{negl}.}
We have
\bqn{\label{in}
\varsigma(u) \leq
\psi_{\delta,H}(u) \leq
\pk{\sup\limits_{t \in G(\delta/u)\backslash I(t_0)}V(t)>u^{1-H}}+
\varsigma(u)
.}
Denote $\ve(t_0) = (-\ve+t_0,\ve+t_0)_{\delta/u}$ for some $\ve>0$.
We have
\bqny{
\pk{\sup\limits_{t \in G(\delta/u)\backslash I(t_0)}V(t)>u^{1-H}} \leq
\pk{\sup\limits_{t \in \ve(t_0)\backslash I(t_0)}V(t)>u^{1-H}}+
\pk{\sup\limits_{t \in G(\delta/u)\backslash \ve(t_0)}V(t)>u^{1-H}}
.}
By Borell-TIS inequality we have
\bqn{\label{bo}
\pk{\sup\limits_{t \in G(\delta/u)\backslash \ve(t_0)}V(t)>u^{1-H}} =
o(\PH(u^{1-H}
\frac{1+ct_0}{t_0^H})), \ \ u \to \IF.
}
Notice, that
\bqny{
\pk{\sup\limits_{t \in \ve(t_0)\backslash I(t_0)}V(t)>u^{1-H}}
&\leq&
Cu\sup\limits_{ t \in \ve(t_0)\backslash I(t_0)}\pk{V(t)>u^{1-H}}
\\&\leq&
Cu\left( \pk{V(t_0-1/\sqrt  u)>u^{1-H}}+
\pk{V(t_0+1/\sqrt  u)>u^{1-H}}\right)
\\&=&
Cu\left(
\PH(u^{1-H}\frac{1+c(t_0-1/\sqrt u)}{(t_0-1/\sqrt u)^H})+
\PH(u^{1-H}\frac{1+c(t_0+1/\sqrt u)}{(t_0+1/\sqrt u)^H})
\right)
\\ &=&
2Cu\PH(u^{1-H}\frac{1+ct_0}{t_0^H})\exp(-\frac{1}{4}f_H''(t_0)u^{1-2H})(1+o(1)), \ \ u \to \IF,
}
recall that $f_H(t) = \frac{(1+ct)^2}{t^{2H}}$ and $f_H''(t_0)>0$. Hence we have
$$\pk{\sup\limits_{t \in \ve(t_0)\backslash I(t_0)}V(t)>u^{1-H}}= o(\PH(u^{1-H}
\frac{1+ct_0}{t_0^H})), \ \ u \to \IF$$
and the claim follows by the line above combined with \eqref{in} and
\eqref{bo}.
\QED
\\
\\
\textbf{Proof of \eqref{bonf}.}
By Bonferroni inequality,
\bqny{
\sum\limits_{t \in I(t_0)}\pk{V(t)>u^{1-H}}-\Pi(u)
\leq
\varsigma(u) \leq \sum\limits_{t \in I(t_0)}\pk{V(t)>u^{1-H}}
,}
where
$$\Pi(u) = \sum\limits_{t_1\neq t_2 \in I(t_0)}\pk{V(t_1)>u^{1-H},
V(t_2)>u^{1-H}}.$$
We show, that
for some positive constant $\hat{C}$ that does not depend on $u$
and for all $t_1,t_2 \in I(t_0)$
\bqn{\label{unif}
\pk{V(t_1)>u^{1-H},
V(t_2)>u^{1-H}} \leq \PH(\hu)\PH(\hat{C}u^{1-2H}), \ \ u \to \IF
}
(recall, $\hu = \frac{u^{1-H}c^H}{H^H(1-H)^{1-H}}$).
Fix some numbers $t_1,t_2 \in I(t_0)$.
Notice that variances of $\frac{V(t_1)}{\sigma_H(t_0)}$ and
$\frac{V(t_2)}{\sigma_H(t_0)}$ are not grater than 1.
Hence
\bqny{\label{upestim}
\pk{V(t_1)>u^{1-H},
V(t_2)>u^{1-H}}
&=&
\pk{\frac{V(t_1)}{\sigma_H(t_0)}>\frac{u^{1-H}}{\sigma_H(t_0)},
\frac{V(t_2)}{\sigma_H(t_0)}>\frac{u^{1-H}}{\sigma_H(t_0)}}
\notag\\ &\leq&
\pk{\frac{\ve^{(1)}V(t_1)}{\sigma_H(t_0)}>\frac{u^{1-H}}{\sigma_H(t_0)},
\frac{\ve^{(2)}V(t_2)}{\sigma_H(t_0)}>\frac{u^{1-H}}{\sigma_H(t_0)}}
\notag \\ &:=&
\pk{W_1>\hu,W_2>\hu},
}
where numbers $\ve^{(1)},\ve^{(2)} \geq 1$ are chosen such that
\bqny{
\Var \frac{\ve^{(1)}V(t_1)}{\sigma_H(t_0)} = \Var
\frac{\ve^{(2)}V(t_2)}{\sigma_H(t_0)} =1
.}
We have, that correlation $r_w$ of $(W_1,W_2)$ has expansion
\bqn{\label{corr}
r_w(t_1,t_2) = 1 - C|t_1-t_2|^{2H}+o(|t_1-t_2|^{2H}), \ \ \ \ t_1,t_2 \to t_0
.}
Notice that size of the grid equals $\delta/u$ and hence for all $t_1,t_2 \in I(t_0)$
\bqn{\label{estim}
\sqrt{|t_1-t_2|^{2H}} \ge \delta^Hu^{-H}  .
}
From Lemma 2.3 in \cite{PickandsA} follows that
\bqny{
\pk{W_1>\hu,
W_2>\hu} \leq
\PH(\hu)\PH(\hu\sqrt{\frac{1-r_w(t_1,t_2)}{1+r_w(t_1,t_2)}}), \ \ u \to \IF.
}
From \eqref{corr} and \eqref{estim} we conclude that for some positive constant $\hat{C}$ as $u \to \IF$
\bqny{
\PH(\hu\sqrt{\frac{1-r_w}{1+r_w}}) \leq \PH(C\hu u^{-H})
=\PH(\hat{C}u^{1-2H}).
}
There are less then $Cu^2$ summands in $\Pi(u)$, hence
by \eqref{unif} and \eqref{fincomput} the claim holds.
\QED
\\
\\
\textbf{Proof of \eqref{comput}.} To establish
\eqref{comput} we need to show that (set $\overline u =
\frac{1}{2}u^{2-2H}$)
\bqny{
\Big|\sum\limits_{t \in I(t_0)}\left(
e^{-\overline u\left(\frac{(1+ct)^2}{t^{2H}} - \frac{(1+ct_0)^2}{t_0^{2H}}\right)}-
e^{-\overline u\frac{f_H''(t_0)}{2}(t-t_0)^2}\right)\Big| = o(1), \ \ u \to \IF.
}
Fix $\rho = \frac{1-2H}{5}>0$. We have
\bqny{ & \ &
\Big|\sum\limits_{t \in I(t_0)}\left(
e^{-\overline u\left(\frac{(1+ct)^2}{t^{2H}} - \frac{(1+ct_0)^2}{t_0^{2H}}\right)}-
e^{-\overline u\frac{f_H''(t_0)}{2}(t-t_0)^2}\right)\Big|
\\&=&
\Big|\sum\limits_{t \in I(t_0)}
e^{-\frac{\overline u}{2}f_H''(t_0)(t-t_0)^2}
\left( e^{-\overline u(\frac{(1+ct)^2}{t^{2H}} - \frac{(1+ct_0)^2}{t_0^{2H}}-
\frac{f_H''(t_0)}{2}(t-t_0)^2)}-1 \right)\Big|
\\ &=&
2\Big|\sum\limits_{t \in (0,u^{-1/2})_{\delta/u}}
e^{-\frac{\overline u}{2}f_H''(t_0)t^2}
\left( e^{-\overline u(\frac{(1+c(t+t_0))^2}{(t+t_0)^{2H}} - \frac{(1+ct_0)^2}{t_0^{2H}}-
\frac{f_H''(t_0)}{2}t^2)}-1 \right)\Big|
\\&=:&
 2\Big|\sum\limits_{t \in (0,u^{-1/2})_{\delta/u}}
e^{-\frac{\overline u}{2}f_H''(t_0)t^2}r(u,t)\Big|
\\ &\leq&
2\Big|\sum\limits_{t \in (0,u^{H-1+\rho})_{\delta/u}}
e^{-\frac{\overline u}{2}f_H''(t_0)t^2}
r(u,t)\Big|
+
2\Big|\sum\limits_{t \in [u^{H-1+\rho},u^{-1/2})_{\delta/u}}
e^{-\frac{\overline u}{2}f_H''(t_0)t^2}
r(u,t)\Big|.
}

We shall show, that both sums above tend to $0$ as $u \to \IF$.
The second sum is negligible because $\overline u t^2 \geq \frac{1}{2}
u^{2\rho}$ for all $t \in [u^{H-1+\rho},u^{-1/2})_{\delta/u}$.
Hence this sum exponentially less then 1 as $u \to \IF$.
For the first sum we write
\bqn{ \label{fs}
\Big|\sum\limits_{t \in (0,u^{H-1+\rho})_{\delta/u}}
e^{-\frac{\overline u}{2}f_H''(t_0)t^2}
r(u,t)\Big|
 \leq
\frac{u^{H+\rho}}{\delta}\sup\limits_{t \in (0,u^{H-1+\rho})}
|r(u,t)|
}
Since $f_H'(t_0) = 0$ and $f_H'''(t_0) \neq 0$
we have as $u \to \IF$
\bqny{
\sup\limits_{t \in (0,u^{H-1+\rho})}
|r(u,t)| \sim \frac{1}{2}u^{2-2H}\frac{|f_H'''(t_0)|}{6}
\sup\limits_{t \in (0,u^{H-1+\rho})}t^3 = Cu^{-1+H+3\rho}.
}
Combining with \eqref{fs} we have as $u \to \IF$
\bqny{
\Big|\sum\limits_{t \in (0,u^{H-1+\rho})_{\delta/u}}
e^{-\frac{\overline u}{2}f_H''(t_0)t^2}
r(u,t)\Big|
 \leq
Cu^{H+\rho-1+H+3\rho}
= Cu^{(-1+2H)/5} \to 0,
}
and the claim is established.\QED

\COM{
the following explicit formula (see \cite{DeM15}):
\bqn{
 \label{nuk}
 \pk{\sup_{t\in [0,T]} ( B(t)- c t)> x}
&=& e^{-2cx}\PH\left(\frac x{ \sqrt{T}} -c\sqrt{T}\right)+
\PH\left( \frac x{ \sqrt{T}} +c\sqrt{T} \right),
}
where $c,T$ and $x$ are any positive numbers.
}

\bibliography{EEEA}{}
\bibliographystyle{plain}
\end{document}